\documentclass[12pt]{amsart}

\newcommand{\BR}{\mathbb{{R}}}

\newtheorem{theorem}{Theorem}[section]
\newtheorem{proposition}[theorem]{Proposition}
\newtheorem{corollary}[theorem]{Corollary}
\theoremstyle{rem}

\theoremstyle{definition}

\def \m{\mu}
\def \l{\lambda}
\def \e{\varepsilon}

\def \o{\omega}
\def \a{\alpha}

\def \d{\delta}
\def \D{\Delta}
\def \g{\gamma}

\def \p{\phi}

\def \b{\beta}

\numberwithin{equation}{section}

\begin{document}

\author{E. Liflyand\qquad and \qquad R. Trigub}
\title {\it On absolute convergence of Fourier integrals}


\address{Department of Mathematics, Bar-Ilan University, 52900 Ramat-Gan, Israel}
\email{liflyand@math.biu.ac.il}
\address{Department of Mathematics, Donetsk National University,
83055 Donetsk, Ukraine}
\email{roald@ukrpost.ua}

\begin{abstract}
New sufficient conditions for representation of a function via the
absolutely convergent Fourier integral are obtained in the paper.
In the main result, Theorem 1.1, this is controlled by the
behavior near infinity of both the function and its derivative.
This result is extended to any dimension $d\ge2.$

\end{abstract}

\maketitle

\section{Introduction}

The possibility to represent a function via the absolutely
convergent Fourier integral was studied by many mathematicians and
is of importance in various problems of analysis. New sufficient
conditions of such type in $\BR^d,$ $d\ge1,$ are obtained in the
paper.

\bigskip

\subsection{Motivation.}

To illustrate the importance of such theorems, let us start with a
typical example of their application. For simplicity, let us
restrict ourselves with the one-dimensional case.

Given algebraic polynomials $P_1,P_2$ and $Q,$ when, for $D=\frac
{d}{dx},$ the inequality

\begin{eqnarray}\label{in1}
\qquad ||Q(D)f||_{L_q(\BR)}\le \g (||P_1(D)f||_{L_{p_1}(\BR)}+
||P_2(D)f||_{L_{p_2}(\BR)})                        \end{eqnarray}
holds true with a constant $\g$ independent of the function $f$?

Let us simplify the situation even more by putting only one
operator on the right. The initial problem is now reduced to when

\begin{eqnarray}\label{in2}
||Q(D)f||_{L_q(\BR)}\le||P(D)f||_{L_p(\BR)}\quad? \end{eqnarray}
Clearly, there should be $s=\mbox{\rm deg}\, Q \le r=\mbox{\rm
deg}\, P$. Considering all the functions for which the right-hand
side of (\ref{in2}) is finite, we deduce that all the solutions of
the equation $P(D)f)=0$ are the solutions of the equation
$Q(D)f=0$ as well, that is, $Q=cP$ with some constant $c.$

The problem becomes meaty under assumption $f\in W^r_p(\BR)$. The
latter is a usual notation for the Sobolev space. Necessary
conditions for fulfillment of (\ref{in2}) in this case are $q\ge
p$ when $s<r,$ $q=p$ when $s=r,$ and when either $q\ne \infty$ or
$p\ne 1$

\begin{eqnarray}\label{in3}
\sup_{x\in \BR}|\p (x)|<\infty,\qquad \p (x)=\frac {Q(ix)}{P(ix)}.
\end{eqnarray}
We have for $|x|\to \infty,$ provided (\ref{in3}) holds true,

\begin{eqnarray*} \p (x)=a_0 + \frac {a_1}{x} +O\Big(\frac {1}{x^2}\Big),\quad \p
^{\prime}(x)= -\frac {a_1}{x^2}+ O\Big(\frac
{1}{x^3}\Big),\end{eqnarray*} and each of the theorems of Section
2 below yields that the Fourier transform $\widehat {\p -a_0} \in
L_1(\BR),$ or in other words $\p -a_0$ is the Fourier transform of
a function $g\in L_1(\BR).$

And if, for example, $f\in W_1^r(\BR),$ then $\widehat
{f^{(k)}}(y)=(iy)^k \widehat f (y),$ $0\le k\le r,$ and

\begin{eqnarray*} \widehat {P(D)f}(y)=P(iy)\widehat f (y),\qquad \widehat
{Q(D)f}(y)=Q(iy)\widehat f (y),\end{eqnarray*}
\begin{eqnarray*}\widehat {Q(D)f}=\p \, \widehat {P(D)f}.\end{eqnarray*}
By this we get that $Q(D)f$ is representable as the convolution of
$g$ and $P(D)f$ (for general theory of multipliers, see \cite{Stein,
SW, St93}; see also the beginning of Section \ref{known}).

In fact, $g$ is also bounded almost everywhere (a.e.). Applying
the Young inequality for convolutions (see, e.g.,
\cite[App.A]{Stein}), we obtain the inequality (\ref{in2}) in the
general case (for details, see \cite{Tr07}, where three criteria
for existence of such inequalities are found: on the axis, on the
half-axis, and on the circle).

We are now in a position to return to the case of two operators
$P_1(D)$ and $P_2(D)$ on the right-hand side of (\ref{in1}).

We suppose that $r_2=\mbox{\rm deg} P_2\le \mbox{\rm deg}
P_1=r_1,$ $P_1(x)=I(x) I_1(x) \widetilde {P_1}(x),$ and
$P_2(x)=I(x) I_2(x) \widetilde {P_2}(x),$ where the polynomials
$\tilde {P_1}$ and $\tilde {P_2}$ never vanish on the imaginary
axis $i\BR$, while zeros of $I, I_1$ and $I_2$, if exist, are
located just on $i\BR.$ Suppose also that $I_1$ and $I_2$ have no
common zeros. In this case the polynomial $Q$ is divisible by $I$
as well.

Obviously, all the values $I_1$ takes on $i\BR$ are on the same
line passing through the origin (as well as those of $I_2 $),
hence one may assume that $I_1$ and $I_2$ take on $i\BR$ only real
value (maybe being multiplied by a constant).

If $I_2(ix_1)\ne 0$ for $x_1\in \BR,$ and $k=\mbox{\rm deg}\,
\widetilde {P_1}$, we set

\begin{eqnarray*}
P_0(x)=\pm iH(x)I_1(x)+I_2(x),\qquad H(x)=(ix+x_1)^k,
\end{eqnarray*} where the sign $+$ or $-$ is chosen in such a way that
$\mbox{\rm deg} P_{0} I=\mbox{\rm deg} P.$ Then $P_0(ix)\ne 0$ for
$x\in \BR.$ Applying (\ref{in2}) three times, we obtain

\begin{eqnarray*} ||Q(D)f||&\le& \g_1||IP_0(D)f||\le
\g_1(||H\,I\,I_1(D)f||+||I\,I_2(D)f||)\\ &\le&
\g_2(||P_1(D)f||+||P_2(D)f||).                   \end{eqnarray*}
We mention that similar arguments are applicable for functions of
several variables as well; more precisely, for elliptic
differential operators and those related. For this, see
\cite{Stein, Bes, BDM} and references therein.

By the way, (\ref{in2}) for $P_1(x)=x^r,$ $r\ge 2,$ $P_2(x)\equiv
1,$ and $Q(x)=x^k,$ $1\le k\le r-1,$ yields, for example for $f\in
W_\infty ^r,$

\begin{eqnarray*}||f^{(k)}||_\infty\le \g_3 (k,r)\,
(||f^{(r)}||_\infty+||f||_\infty ).     \end{eqnarray*} Substituting
$\e x,$ with $\e >0,$ for $x,$ dividing by $\e^k$ and minimizing the
right-hand side over $\e\in (0,\infty),$ we obtain the known
multiplier inequality for intermediate derivatives

\begin{eqnarray*}
||f^{(k)}||_\infty\le \g _4 (k,r)||f||^{\frac kr}_\infty\,
||f^{(r)}||^{1-\frac kr} _{\infty}. \end{eqnarray*}

Prior to formulating main results we explain how the paper is
organized and fix certain notation and conventions. First, if

\begin{eqnarray*}
f(y)=\int\limits_{\BR^d}g(x)e^{i(x,y)}dx,\qquad g\in L_1(\BR^d),
\end{eqnarray*}
we write $f\in A(\BR^d),$ with $||f||_A=||g||_{L_1(\BR^d)}.$ In
Section \ref{known}, known results on representability of a function
as the absolutely convergent Fourier integral are given and
compared. The proofs of the new results are given in Section
\ref{proofs}.

We shall denote absolute constants by $c$ or maybe by $c$ with
various subscripts, like $c_1,$ $c_2,$ etc., while $\g(...)$ will
denote positive quantities depending only on the arguments indicated
in the parentheses. We shall also use the notation $\int_{\to0}$ to
indicate that the integral is understood as improper in a
neighborhood of the origin, that is, as
$\lim\limits_{\delta\to0+}\int_\delta.$

\bigskip

\subsection{Main results.}

We start with the case $d=1$.

Let $f\in C_0(\BR),$ that is, $f\in C(\BR)$ and $\lim f(t)=0$ as
$|t|\to \infty$, and let $f$ be locally absolutely continuous on
$\BR\setminus\{0\}.$

\begin{theorem}\label{th1}
Let $f_0(t)=\sup_{|s|\ge |t|}|f(s)|.$

{\bf a)} Let $f^{\prime}$ be bounded a.e. out of any neighborhood of
zero and $f_1(t)=\mathrel{\mathop{\mbox{\rm ess\,sup}}}_{|s|\ge
|t|>0} |f^{\prime}(s)|$. If, in addition,

\begin{eqnarray*}
A_0=\int\limits_1^\infty \frac {f_0(t)}{t}\,dt<\infty,\qquad
A_1=\int\limits_0^1 f_1(t)\ln\frac 2t\, dt <\infty
\end{eqnarray*}
and

\begin{eqnarray*} A_{01}=\int\limits_1^\infty \biggl(\int\limits_t^\infty
f_0(s)f_1(s)\,ds\biggr)^\frac12\,\frac {dt}{t}<\infty,
\end{eqnarray*} then $f\in A(\BR),$ with $||f||_A\le c
(A_0+A_1+A_{01}).$

{\bf b)} Let  $f^\prime$ be not bounded near infinity, $f_\infty(t)=
\mathrel{\mathop{\mbox{\rm ess\,sup}}}_{0<|s|\le |t|}|f^\prime (s)|$
and $f(t)=0$ when $|t|\le 2\pi,$ with $f_\infty(4\pi)>0.$ If, in
addition, there exists $\d\in (0,1)$ such that

\begin{eqnarray*} A_{\d}^{1+\d}=\sup_{t\ge 2\pi} t\, f_0^\d(t)\,
f_\infty(t+2\pi)<\infty,                       \end{eqnarray*} then
$f\in A(\BR)$ and $||f||_A \le \g(\d)A_\d (1+A_{\d}^{\frac
1\d}(f_\infty(4\pi))^{-\frac 1\d}).$
\end{theorem}

Conditions of this theorem differ from known sufficient conditions
in the way that near infinity combined behavior of both the function
and its derivative comes into play (see also the corollary below).
Conditions for $f_0$ near infinity and $f_1$ near the origin in {\bf
a)} are also necessary. For instance, this is the case when $f(t)=0$
for $t\le 0,$ $f\in C^1(0,+\infty)$  and piecewise convex on
$[0,\infty),$ since for such functions both conditions are
equivalent to convergence of the integral $\int\limits_0^\infty
t^{-1}f(t)\,dt$ (see necessary conditions in Section 2). As for the
condition $A_{01}<\infty,$ it holds, for example, if $(\ln
t)^{2+\delta} f_0(t)f_1(t) \in L_1[1,\infty)$ for some $\delta
>0,$ but nit for $\delta=0.$

We note that in {\bf b)} the function can be considered on the
whole axis. It should satisfy the same condition $A_1$ as in {\bf
a)} near the origin. We omit this for simplicity.

\begin{corollary}\label{sl} If $A_1<\infty,$ $f(t)=O(|t|^{-\a})$
for some $\a>0$ and $f^\prime (t)=O(|t|^{-\b})$ for some $\b\in \BR$
as $|t|\to \infty,$ with $\a+\b>1,$ then $f \in A(\BR).$ If $\a+\b
<1$ such an assertion cannot be valid.
\end{corollary}

Let now $f:\BR^d\to \mathbb C$ with $d\ge 2.$ We will give a
direct generalization of {\bf a)} in Theorem \ref{th1} to higher
dimensions.

To formulate a multivariate extension of that result, we introduce
certain notation. Let $\chi,$ $\eta$ and $\zeta$ be $d$-dimensional
vectors with the entries either $0$ or $1$ only. Each of these
vectors or even two of them can be zero vectors ${\bf
0}=(0,0,...,0).$ The inequality of vectors means the same inequality
for all pairs of their corresponding components.

Similarly to the one-dimensional case, we set for $\eta+\zeta={\bf
1}=\{1,1,...,1\}.$

\begin{eqnarray*} f_{\eta,\zeta}(x)=\sup\limits_{|u_i|\ge|x_i|,\atop i:\eta_i=1}
\mathrel{\mathop{\mbox{\rm ess\,sup}}_{|u_j|\ge |x_j|>0,\atop
j:\zeta_j=1}}|D^{\zeta} f(u)|,         \end{eqnarray*} where

\begin{eqnarray*} D^\chi f(x)=\left(\prod\limits_{j: \chi_j=1}
\frac{\partial}{\partial x_j}\right)f(x).\end{eqnarray*}

We denote by $\mathbb R_\zeta$ the Euclidean space of dimension
$\zeta_1+...+\zeta_d$ with respect to the variables $x_j$ with
$j$s for which $\zeta_j=1;$ correspondingly $x_\zeta$ is an
element of this space.

\begin{theorem}\label{newmu}
Let $f\in C_0(\mathbb R^d)$ and let $f$ and its partial derivatives
$D^\eta f,$ ${\bf 0}\le\eta<{\bf 1},$ be locally absolutely
continuous on $(\mathbb R\setminus \{0\})^d$ in each variable. Let
also partial derivatives $D^\eta f,$ ${\bf 0}<\eta\le{\bf 1}$ be
almost everywhere bounded out of any neighborhood of each coordinate
hyperplane. If

\begin{eqnarray}\label{muco}&\quad&A_{\chi,\eta,\zeta}=\int\limits_0^1...\int\limits_0^1
\prod\limits_{k:\zeta_k=1}\ln(2/x_k)\,dx_k\,\nonumber\\
&\quad&\int\limits_1^\infty...\int\limits_1^\infty
\biggl(\int\limits_{\prod\limits_{j:\eta_j=1}[u_j,\infty)}f_{\chi+\eta,\zeta}(x)
f_{\chi,\eta+\zeta}(x)\,dx_\eta\biggr)^{1/2}\prod\limits_{i:\chi_i=1\atop
{\rm or}\ \eta_i=1}\frac {du_i}{u_i}<\infty       \end{eqnarray}
for all $\chi,$ $\eta$ and $\zeta$ such that $\chi+\eta+\zeta={\bf
1}=\{1,1,...,1\},$ then $f\in A(\mathbb R^d).$ \end{theorem}

The authors understand, of course, that (\ref{muco}) is a (quite
large) number of conditions not easily observable, in a sense. To
clarify this point, we give the two-dimensional version of this
theorem. We will use separate letter for each variable rather than
subscripts; also no need in using vectors for defining analogs of
$f_0$ and $f_1$ - we just denote them by using subscripts $0$ or
$1$ to indicate majorizing in the corresponding variable:
$f_{00},$ $f_{01},$ $f_{10}$ and $f_{11}.$

{\bf Theorem \ref{newmu}$'$.} {\it Let $f(x,y)\in C_0(\mathbb R^2)$
and let $f$ and its partial derivatives $\frac{\partial f}{\partial
x}$ and $\frac{\partial f}{\partial y}$ be locally absolutely
continuous on $(\mathbb R\setminus \{0\})^2$ in each variable. Let
also partial derivatives $\frac{\partial f}{\partial x},$
$\frac{\partial f}{\partial y}$ and $\frac{\partial^2 f}{\partial
x\partial y}$ be almost everywhere bounded out of any neighborhood
of each coordinate axis. If

\begin{eqnarray*}\int\limits_0^1\int\limits_0^1 f_{11}(x,y)
\ln\frac{2}{x}\ln\frac{2}{x}\,dx\,dy<\infty,\quad
\int\limits_1^\infty\int\limits_1^\infty
\frac{f_{00}(x,y)}{xy}\,dx\,dy<\infty,\end{eqnarray*}

\begin{eqnarray*}\int\limits_1^\infty\int\limits_1^\infty\biggl(
\int\limits_x^\infty\int\limits_y^\infty f_{00}(s,t)f_{11}(s,t)
\,ds\,dt\biggr)^{1/2}\frac{dx}{x}\frac{dy}{y}<\infty,\end{eqnarray*}

\begin{eqnarray*}\int\limits_0^1\int\limits_1^\infty f_{10}(x,y)
\ln\frac{2}{x}\,dx\,\frac{dy}{y}<\infty,\end{eqnarray*}

\begin{eqnarray*}\int\limits_0^1\int\limits_1^\infty \biggl(\int\limits_y^\infty
f_{10}(x,t)f_{11}(x,t)\,dt\biggr)^{1/2}
\ln\frac{2}{x}\,dx\,\frac{dy}{y}<\infty,\end{eqnarray*}

\begin{eqnarray*}\int\limits_1^\infty\int\limits_0^1 f_{01}(x,y)
\ln\frac{2}{y}\,\frac{dx}{x}\,dy<\infty,\end{eqnarray*}

\begin{eqnarray*}\int\limits_1^\infty\int\limits_0^1 \biggl(\int\limits_x^\infty
f_{01}(s,y)f_{11}(s,y)\,ds\biggr)^{1/2}
\ln\frac{2}{y}\,\frac{dx}{x}\,dy<\infty,\end{eqnarray*}

\begin{eqnarray*}\int\limits_1^\infty\int\limits_1^\infty\biggl(
\int\limits_y^\infty f_{00}(x,t)f_{01}(x,t)\,dt\biggr)^{1/2}
\,\frac{dx}{x}\,\frac{dy}{y}<\infty,\end{eqnarray*} and

\begin{eqnarray*}\int\limits_1^\infty\int\limits_1^\infty\biggl(
\int\limits_x^\infty f_{00}(s,y)f_{10}(s,y)\,ds\biggr)^{1/2}
\,\frac{dx}{x}\,\frac{dy}{y}<\infty,           \end{eqnarray*} then}
$f\in A(\mathbb R^2).$

As is mentioned, already for $d=3$ no way to briefly write down
all the conditions. Let is give only one of them, quite typical
and completely "mixed":

\begin{eqnarray*}\int\limits_0^1\int\limits_1^\infty\int\limits_1^\infty\biggl(
\int\limits_z^\infty f_{100}(x,y,u)f_{101}(x,y,u)\,du\biggr)^{1/2}
\,dx\,\frac{dy}{y}\,\frac{dz}{z}<\infty;        \end{eqnarray*}
here $f_{100}$ and $f_{101}$ is a clear analog of the above
notation.

\bigskip

\section{Known results}\label{known}

Let $\p:\BR^d \to \mathbb C$ be a bounded measurable function. We
define on $L_2(\BR^d)\cap L_p(\BR^d)$ a linear operator $\Phi$ via
the following equality for the Fourier transform of a function
$f\in L_2\cap L_p$

\begin{eqnarray*}
\widehat {\Phi f}(y)=\p(y)\widehat{f}(y).  \end{eqnarray*}
Clearly, $\Phi f\in L_2 (\BR^d),$ and if there exists a constant
$D$ such that for all $f\in L_2 \cap L_p (\BR^d)$

\begin{eqnarray*}||\Phi (f)||_{L_p}\le D ||f||_{L_p}, \end{eqnarray*}
then the operator $\Phi$ is called the Fourier multiplier from
$L_p$ into $L_p(\BR^d)$ (written $\p\in M_p(\BR^d)$), whence
$||\Phi||_{L_p\to L_p}=\inf D.$

Sufficient conditions for a function to be a multiplier in $L_p$
spaces with $1<p<\infty$ for both multiple Fourier series and
Fourier integrals were studied by Marcinkiewicz, Mikhlin,
H\"ormander, Lizorkin, and others (see, e.g., \cite[Ch.4]{Stein} and
\cite{St93}). There holds $M_1=M_\infty \subset M_p,$ $1<p<\infty.$
When $p=1$ and $p=\infty$ each Fourier multiplier is the convolution
of the function $f$ and a finite (complex-valued) Borel measure on
$\BR^d$:

\begin{eqnarray*}
\Phi f(x)=\int\limits_{\BR^d}f(x-y)d\m(y),\qquad ||\Phi||_{L_1\to
L_1}=||\Phi||_{L_\infty \to L_\infty}=\mbox{\rm var} \m,
\end{eqnarray*}
while $\p\in M_1=M_\infty$ iff $\p\in B(\BR^d),$ where

\begin{eqnarray*}
B(\BR^d)=\{\p:\p(y)=\int\limits_{\BR^d}e^{i(x,y)}d\m (x),\qquad
||\p||_B=\mbox{\rm var} \m<\infty\}
\end{eqnarray*}
(see, e.g., \cite[Ch.1]{SW}). If the measure $\m$ is absolutely
continuous with respect to the Lebesgue measure in $\BR^d,$ then
we write $\p\in A(\BR^d),$ where

\begin{eqnarray*}
A(\BR^d)&=&\{\p :\p (y)=(2\pi)^{\frac d2}\widehat g
(-y)=\int\limits_{\BR^d}e^{i(x,y)}g(x)dx, \\ ||\p||_A
&=&\int\limits_{\BR^d}|g(x)|dx<\infty \}.
\end{eqnarray*} The space $B(\BR^d)$ is the Banach algebra with
respect to pointwise multiplication, while $A(\BR^d)$ is an ideal
in $B(\BR^d).$ As is known, the two algebras are locally geared in
the same way, thus the difference between $A$ and $B$ is revealed
in the behavior of functions near infinity. We also note that if
$\p\in B(\BR^d),$ $\lim \p (y)=0 $ as $|x|\to \infty$ and $\p$ is
of finite total Vitali variation off a cube, then $\p \in
A(\BR^d)$ (\cite[Theorem 2]{Tiz}).

We remind the reader that total Vitali variation of the function
$\p:E\to \mathbb C,$ with $E\ \subset \BR^d,$ is defined as
follows. If $\{e_j^0\}_{j=1}^d$ is the standard basis in $\BR^d$,
and the boundary of $E$ consists of a finite number of planes
given by equations $x_j=c_j,$ then

\begin{eqnarray*} V(f)=\sup\sum|\D_uf(x)|,\qquad
\D_uf(x)=(\prod\limits_{j=1}^{d}\D_{u_j})f(x),
\end{eqnarray*}
where $u=(u_1,\cdots ,u_d)$ and

\begin{eqnarray}\label{in4}
\D_{u_j}f(x)=f(x+u_j e_j^0)-f(x-u_j e_j^0), \quad 1\le j\le d.
\end{eqnarray}
Here $\D_u$ is the mixed difference with respect to the vertices
of the parallelepiped $[x-u,x+u]$ and $\sup$ is taken over any
number of non-overlapping parallelepipeds in $E.$ For smooth
enough functions on $E$ such as indicated above, one has

\begin{eqnarray*}
V(f)=\int\limits_E \left|\frac {\partial ^d f(x)}{\partial x_1
\cdots\partial x_d}\right|\,dx.                  \end{eqnarray*}
We note that in Marcinkiewicz's sufficient condition for $M_p,$
$1<p<\infty,$ only the finiteness of total variations over all
dyadic parallelepipeds with no intersections with coordinate
hyper-planes is assumed (see, e.g., \cite{Stein}).

Many mathematicians studied the properties of absolutely
convergent Fourier series rather than integrals, starting from one
paper by S.N. Bernstein (see, e.g., \cite{Kahane}; for
multidimensional results see, e.g., \cite{Tim}). Various
sufficient conditions for absolute convergence of Fourier
integrals were obtained by Titchmarsh, Beurling, Karleman,
Sz.-Nagy, Stein, and many others. One can find more or less
comprehensive and very useful survey on this problem in \cite{SK},
with 65 bibliographical references therein.

Let us give some results not contained in that survey as well as
relations between them and other results of such type. The other
reason for giving these is that some of these results will
essentially be used in proofs.

P\'olya proved that each even, convex and monotone decreasing to
zero function on $[0,\infty)$ belongs to $A(\BR).$ In fact, such
function belongs even to $A^*(\BR),$ that is, not only $\widehat
f\in L_1(\BR),$ but also $\sup_{|s|\ge |t|}|\widehat f (s)|\in
L_1(\BR)$ (see \cite{Tiz} or \cite{TB}). By this, $f$ may decrease
arbitrarily slowly. What is really important, as P\'olya observed,
is that $\widehat f (y)\ge 0.$

Zygmund  proved that if an odd function $f$ is compactly supported
and convex in a right neighborhood of the origin, it admits an
extension to $A(\BR)$ iff the improper integral $\int_{\to0}
t^{-1}f(t)\,dt $ converges (see \cite{Kahane}). There is a more
general statement (Lemma 6 in \cite{Tiz}): if $f\in C_0(\BR)$ and
piece-wise convex, then for any $y\ne 0$

\begin{eqnarray*}|\widehat f (y)|\le\g(q)\o(f;\frac {\pi}{|y|}),
\end{eqnarray*}
where $\o(f;h)$ is the modulus of continuity and $q$ is the number
of intervals on each of them $\mbox{\rm Re} f$ and $\mbox{\rm Im}
f$ are either convex or concave. It follows from this that if a
function is also odd and on any interval not containing zero
satisfies the $\mbox{\rm Lip} \a,$ $\a>0,$ condition, then $f\in
A(\BR)$ iff the integral $\int_0^\infty t^{-1}f(t)\,dt$ converges.
Therefore, if $f$ is odd and $f(x)\ge 0$ for $x\ge 0,$ then $f\in
A(\BR)$ does not yield $f\in L_1(\BR)$ (cf. Theorem 2.8 in
\cite{SK}).

For a real, bounded and locally absolutely continuous function to be
the difference of two bounded convex functions on $[0,\infty),$
(quasi-convex), it is necessary and sufficient that

\begin{eqnarray*}\int\limits_0^\infty t|df^{\prime}(t)|<\infty.\end{eqnarray*}
A similar fact is well known for sequences.

In the paper by Beurling \cite{Beur} more general condition was
given:

\begin{eqnarray*}
V^*(f)=\int\limits_0^\infty \mathrel{\mathop{\mbox{\rm
ess\,sup}}_{s\ge t}} |f^{\prime}(s)|\,dt<\infty.
\end{eqnarray*}
This condition is less restrictive that that of convexity(and
quasi-convex), but more severe than the finiteness of the total
variation.

If, in addition, $f\in C_0[0,\infty)$ and $f(t)=0$ for $t<0,$ then
for each $y\in \BR\ \setminus \{0\}$

\begin{eqnarray*}\widehat f(y)&=&\frac {1}{\sqrt {2\pi}}\int _0^\infty
f(t)e^{-iyt}dt=-\frac{i}{y\sqrt{2\pi}}f(\frac{\pi}{2|y|})+\theta
F(y),
\end{eqnarray*}
with $|\theta|\le c$ and $||F||_{L_1(\BR)}\le V^*(f).$

This can be found in \cite[6.4.7b and 6.5.9]{TB}. It is
interesting that for $f$ convex $F(y)$ can be considered monotone
decreasing as $|y|$ increases, while it is not the case for the
class in question. The point is that integrability of the monotone
majorant of $|\widehat f|$ is related to summability of Fourier
series at Lebesgue points (see \cite[8.1.3]{TB}). Monotonicity of
$F$ in the latter case will then lead to coincidence of results
for $V^*$ and convex functions, which is impossible.

{\bf Theorem A} {\rm (\cite{Beur})}. Let $f\in C_0(\BR)$ and there
exists a function $g$ such that

\begin{eqnarray*}
|f(t)-f(t+h)|\le |g(t)-g(t+h)| \qquad (t,h\in \BR)
\end{eqnarray*}
and $g\in A^*(\BR)$, that is, $g=\widehat \psi,$ with
$\psi^*(t)=\mathrel{\mathop{\mbox{\rm ess\,sup}}_{|t|\ge
|s|}}|\psi (s)|\in L(\BR).$ Then $f\in A(\BR).$

For general properties of the algebra $A^*(\BR),$ see \cite{BLT}.

On the other hand, many works were devoted to the related question
of the boundedness and asymptotics of $L_1$-norms over the period
of the sequence of periodic functions against their Fourier
coefficients (see, e.g., \cite{Te}, and also  \cite[7.2.8,
8.1.1]{TB} and \cite{L2}).

Let us now give F. Riesz's criterion of the absolute convergence
of Fourier integrals (its counterpart for series can be found in
\cite{Kahane}).

{\bf Theorem B.} {\it Function $f\in A(\mathbb R^d)$ if and only
if it is representable as the convolution of two functions from
$L_2 (\mathbb R^d)$:

\begin{eqnarray}\label{conv}
f(x)=\int\limits_{\BR^d}f_1(y)f_2(x-y)dy,\qquad f_1,f_2\in
L_2(\BR^d).
\end{eqnarray}
By this}, $||f||_A\le ||f_1||_2\,||f_2||_2.$

\begin{proof} The proof is based on the unitarity of the Fourier
operator in $L_2(\BR^d)$. With (\ref{conv}) in hand,

\begin{eqnarray*}
f(x)=\int\limits_{\BR^d}\widehat {f_1}(u)\overline{\widehat
{f_2(x-u)}}du=\int\limits_{\BR^d} \widehat {f_1}(u)\, \overline
{\widehat {f_2}(-u)}\, e^{i(x,u)}du
\end{eqnarray*} and by the Cauchy-Schwarz-Bunyakovskii inequality

\begin{eqnarray*}||f||_A=\int\limits_{\BR^d}|\widehat {f_1}(u)|\,|\widehat
{f_2}(-u)|\,du\le||\widehat{f_1}||_2\,||\widehat{f_2}||_2=||f_1||_2\,||f_2||_2.
\end{eqnarray*}
Farther, if

\begin{eqnarray*} f(x)=\int\limits_{\BR^d}g(y)e^{i(x,y)}dy,\qquad
||g||_1=\int\limits_{\BR^d}|g(y)|\,dy<\infty,     \end{eqnarray*}
then

\begin{eqnarray*}   g(y)=|g(y)|^{\frac 12}
\,(|g(y)|^{\frac 12}\,\mbox{\rm sign}\,g(y)),\end{eqnarray*}
where, as usual, $\mbox{\rm sign} 0=0$ and $\mbox{\rm sign}z=
\frac{z}{|z|}$ when $z\ne 0,$ and each of the two factors on the
right, as well as their Fourier transforms, belongs to
$L_2(\BR^d).$ It remains to apply the above given formulae for
convolution in the reverse order. \hfill\end{proof}

Let us demonstrate how to derive effective sufficient conditions
from this criterion.

Assume $f\in C(\BR^d)\cap L_2(\BR^d)$ and $(-\D)^{\frac \a 2}f\in
L_2(\BR^d),$ where $\D$ is the Laplace operator. Taking into account
that

\begin{eqnarray*}
\widehat{(-\D)^{\alpha/2}f}(y)=|y|^\a\widehat{f}(y)\in L_2(\BR^d),
\end{eqnarray*}
and for $\a>d/2$

\begin{eqnarray*}
\int\limits_{\BR^d}\frac
{dy}{(1+|y|)^2}=\g(\a)\int\limits_0^\infty \frac {t^{\a-1}}
{(1+t^\a)^2}dt<\infty,
\end{eqnarray*}
we get the product of the two functions from $L_{2} (\BR^d)$

\begin{eqnarray*}
\widehat{f}(y)=(\widehat{f}(y)(1+|y|^\a))\, \frac {1}{1+|y|^\a}.
\end{eqnarray*}
Therefore $f\in A(\BR^d)$ when $\a> d/2$.

Another differential operators can be used in the same way, say
elliptic, while applying embedding theorems allows one to digress
on the function classes defined via moduli of continuity of
partial derivatives.

There is one more criterion (approximative) from which in \cite
[6.4.3]{TB}, for example, known sufficient conditions with
different smoothness in various variables are derived.

Let us go on to necessary conditions for dimension one. Obviously,
for each $\l\in \BR$

\begin{eqnarray*}
||\overline f||_A=||f(\l\cdot)||_A=||f(\cdot+\l)||_A=
||e^{i\l(\cdot)}f(\cdot)||_A =||f||_A.
\end{eqnarray*}
let $f\in A(\BR)$ and $f=\sqrt{2\pi}\widehat g,$ where $g\in
L_1(\BR)$. Then the trigonometrically conjugate function (the
Hilbert transform) is

\begin{eqnarray*}\widetilde f(x)&=&\frac 1\pi\int\limits_{\to0}^{\to \infty} \frac
{f(x+t)-f(x-t)}{t}dt\\
&:=& \lim_{\e\to+0,M\to +\infty}\frac 1\pi \int\limits_\e^M \frac {f(x+t)-f(x-t)}{t}dt\\
&=& \lim_{\e\to+0,M\to +\infty}\frac 1\pi \int\limits_\e^M
dt\int\limits_{-\infty}^{+\infty}g(y)\frac{e^{iy(x+t)}-e^{iy(x-t}} {t}dy\\
&=&\frac {2i}{\pi}\lim_{\e\to+0,M\to +\infty}
\int\limits_{-\infty}^{+\infty}g(y)e^{ixy}dy \int\limits_\e^M\frac
{\sin ty}{t}dt.                \end{eqnarray*} Since the absolute
values of the integrals over $[\e,M]$ are bounded by an absolute
constant, it is possible to pass to the limit under integral sign.
This yields

\begin{eqnarray*}
\widetilde
{f}(x)=i\int\limits_{-\infty}^{+\infty}g(y)e^{ixy}\mbox{\rm sign}
y dy,\qquad ||\widetilde f||_A=||f||_A.
\end{eqnarray*}
We mention that the improper integral in the definition of
$\widetilde f$ converges everywhere (and uniformly in $x$), but
not necessarily absolutely.

In \cite[Theorem 3]{Tiz}, a necessary condition for belonging to
$A(\BR^d)$ is given. It is valid for both radial and non-radial
functions of $d$ variables and depends on dimension $d.$

To formulate the next result on which much in the proofs of our
new results is based on, we remind that $\D_u f=\D_{u_1,\cdots
,u_d}f $  is defined by (\ref{in4}).

{\bf Theorem C} (Lemma 4 in \cite{Tiz}). {\it Let $f\in
C_0(\BR^d)$.

{\bf a)} If

\begin{eqnarray*}
\sum_{s_1=-\infty}^\infty \cdots \sum_{s_d=-\infty}^\infty 2^ {\frac
12 \sum_{j=1}^d s_j}||\D_{\frac{\pi}{2^{s_1}},\cdots,
\frac{\pi}{2^{s_d}}}(f)||_2<\infty,
\end{eqnarray*}
where the norm is that in $L_2 (\BR^d),$ then $f\in A(\BR^d)$.

{\bf b)} If $f=\widehat g,$ with $g\in L(\BR^d),$ and for
$|u_j|\ge |v_j|$ when $\mbox{\rm sign} u_j=\mbox{\rm sign} v_j$
for all $1\le j \le d$ there holds $|g(u)|\le |g(v)|,$ then the
series in {\bf a)} converges.}

Let us comment on this result.

 The convergence of the series condition in {\bf a)} is of
 Bernstein type. The second named author has learned from \cite{Bes}
 that assertions very similar to {\bf a)} of Theorem C were earlier
 obtained in \cite{Gab} and \cite{Gol} (these results are also
 commented in \cite{Bes}).Besides that, the same assertion {\bf a)}
 is reproved in \cite[Theorem 3]{Bes}; the author had apparently been
 unaware of existence of \cite{Tiz}. Farther, more general problem
 with power weights is studied in \cite{Liz4}.

 However, it is worth mentioning that {\bf b)} is proved in \cite{Tiz}
 as well. Not only this is a necessary condition for certain subclass,
 but application of this assertion to extension of Beirling's theorem
 (see Theorem A above) on functions of any number of variables is given
 in the same paper \cite{Tiz}. To this extent, mixed differences are
 used instead of those simple. In turn, the following sufficient
 condition for belonging to $A(\BR^d)$ was derived from that generalized
 Beurling theorem.

{\bf Theorem D} \mbox{\rm (\cite{Tiz} or \cite{TB})}. {\it If
$f(x)=\int\limits_{|x_j|\le |u_j|,\atop 1\le j\le d} g(u)du,$ with

\begin{eqnarray*}
\int\limits_{\BR^d} \mathrel{\mathop{\mbox{\rm
ess\,sup}}_{|u_j|\ge |v_j|,1\le j\le d}}|g(u)|\,dv<\infty,
\end{eqnarray*} then} $f\in A(\BR^d)$.

As is known, the Hausdorff-Young inequality gives sufficient
condition for $\widehat f \in L_p,$ $p>2.$ In \cite[6.4.2]{TB}
both {\bf a)} and {\bf b)} of Theorem C are generalized to the
case where $\widehat f \in L_1\cap L_p,$ $p\in (0,2).$ We also
mention that when $f(x)=0$ for $x \in\BR\setminus [0,\pi]$ and
$\widehat f \in L_1\cap L_p(\BR),$ $p\in (0,1]$ the following
necessary condition holds true (see the same reference):

\begin{eqnarray*} \int\limits_0^\pi(|f(x)|^p+|f(\pi-x)|^p) x^{p-2}dx
\le \g (p)||\widehat  f||_p^p.                  \end{eqnarray*}
Let us give a recent simple sufficient condition.

{\bf Theorem E} \mbox{\rm (\cite[Theorem 1]{BDM}).}

If $f\in C(\BR^d)$, for any $\d=(\d_1,\cdots ,\d_d),$ with
$\d_j=0$ or $1,$ $1\le j\le d,$

\begin{eqnarray*}
\lim_{|x_j|\to \infty} \frac {\partial ^{\sum \d_j}f(x)}{\partial
x_1^{\d_1} \cdots \partial x_d^{\d_d}}=0, \qquad 1\le j \le d,
\end{eqnarray*}
and for some $\e\in (0,1)$

\begin{eqnarray*}
\bigg|\frac {\partial ^d f(x)}{\partial x_1\cdots \partial x_d} \bigg|\le
\frac {A}{\Pi_{j=1}^d |x_j|^{1-\e}(1+|x_j|)^{2\e}},
\end{eqnarray*}
then $||f||_{A} \le A\g(\e).$

This theorem is a simple consequence of Theorem D. To make the
paper self-contained, let us present this argument.

\begin{proof}[Proof of Theorem E]
We first mention that integrating over $[x_j,\infty)$ an
inequality for the mixed derivative, one obtains a similar
inequality for derivatives of smaller order.

Let, for simplicity, $d=2$. If $f$ as a function of $x_1$ i $x_2$
is even in both $x_1$ and $x_2,$ then

\begin{eqnarray*}
f(x_1,x_2)=\int _{|x_1|}^\infty du\int _{|x_2|}^\infty \frac
{\partial ^2 f(u,v)} {\partial u \partial v}\,dv,
\end{eqnarray*}
and Theorem D is immediately applicable.

Any function $f$ is representable as a sum of at most 4 summands,
each of them is a function either even or odd in $x_1$ and $x_2$
that satisfies the assumptions of Theorem E.

Let, for example, $f_1(-x_1,x_2)=-f_1(x_1,x_2)$, while
$f_1(x_1,-x_2)=f_1(x_1,x_2).$ The function

\begin{eqnarray*}
h_\mu(t)=\begin{cases}|t|^\mu \mbox{\rm sign} t, & -1\le t \le
1,\\|t|^{-\mu} \mbox{\rm sign} t, & |t|>1,\end{cases}
\end{eqnarray*}
belongs to $A(\BR)$ for any $\m>0$. Clearly, always $A(\BR)
\subset B(\BR^2)$.The function $\frac {f}{h_\mu}$ also satisfies
assumptions of Theorem E for $\mu>0$ small enough. As a function
even in both $x_1$ and $x_2$ it belongs to $A(\BR^2)$ due to above
argument. Hence

\begin{eqnarray*}
f_1=\Big(\frac {f_1}{h_\m}\Big)\, h_\m\in A(\BR^2),
\end{eqnarray*}
which completes the proof of Theorem E. \hfill\end{proof}

We mention that there exist results similar to Theorem D for
even-odd and just odd functions (general case) in \cite{GM1}.

But even more general asymptotic formulae for the Fourier
transform can be found in \cite{L0}.

In the problems of integrability of the Fourier transform the
following $T$-transform of a function $h(u)$ defined on
$(0,\infty)$ is of importance

\begin{eqnarray}\label{Ttr} Th(t)=\int\limits_{\to0}^{t/2}
\frac{h(t+s)-h(t-s)}{s}\,ds.                        \end{eqnarray}
In \cite{Fr} it is called the Telyakovskii transform. The reason is
that in an important asymptotic result for the Fourier transform
\cite{L0} (cf. \cite{Fr}) it is used to generalize Telyakovskii's
result for trigonometric series.

It is clear that the $T$-transform should be related to the
Hilbert transform; this is revealed and discussed in \cite{L0} and
later on in, e.g., \cite{Fr}, \cite{L2}, etc. In particular, the
space that proved to be of importance is, related to the real
Hardy space, that of functions with both derivative and
$T$-transform of the derivative being integrable.

We need additional notation different from that in \cite{L0} and
better, in our opinion. Let us denote by $T_jg(x)$ the
$T$-transform of a function $g$ of multivariate argument with
respect to the $j$-th (single) variable:

\begin{eqnarray*}T_jg(x)=\int\limits_{x_j/2}^{3x_j/2}\frac{g(x)}{x_j-t}\,dt
=\int\limits_0^{x_j/2}\frac{g(x-te_j^0)-g(x+te_j^0)}{t}\,dt.\end{eqnarray*}

Analyzing the proof of Theorem 8 in \cite{L0}, one can see that
this theorem can be written in the following asymptotic form.

{\bf Theorem F.} {\it Let $f$ be defined on $\mathbb{R}^d_+;$ let
all partial derivatives taken one time with respect to each
variable involved be locally absolutely continuous with respect to
any other variable all of these derivatives vanish at infinity as
$x_1+...+x_d\to\infty.$ Then for any $x_1,...,x_d>0$ and for any
set of numbers $\{a_j: a_j=0\ \text{or}\ 1\}$ we have

\begin{eqnarray}\label{gam} &\quad&\int\limits_{\mathbb{R}^d_+}f(u)
\prod\limits_{j=1}^d\cos(x_ju_j-\pi a_j/2)\,du_j \nonumber\\
&=&(-1)^{d-1}f\Big(\frac{\pi}{2x_1},...,\frac{\pi}{2x_d}\Big)
\prod\limits_{j=1}^d\frac{\sin(\pi a_j/2)}{x_j}\nonumber\\
&+&\sum\limits_{\chi+\eta+\zeta={\bf1}, \atop \chi\ne{\bf 0}}\,
\int\limits_{\mathbb{R}^d_+}\prod\limits_{i:\chi_i\ne0}\frac{\sin(\pi
a_j/2)}{x_j}F_{\eta+\zeta}(x)\,dx,                  \end{eqnarray}
where $F_{\eta+\zeta}$ functions satisfying}

\begin{eqnarray}\label{esremd}&\quad&\int\limits_{\mathbb{R}^d_+}
|F_{\eta+\zeta}(x)|\,dx\nonumber \\
&\le&c\int\limits_{\mathbb{R}^d_+}\prod\limits_{i:\chi_i\ne0}\frac{\sin(\pi
a_j/2)}{x_j}\biggr|\prod\limits_{j:\eta_j\ne0} T_j\prod
\limits_{k:\tau_k\ne0}D_{k}f(x)\biggr|\,dx.\end{eqnarray}

Integrating the summands in this theorem, we obtain

\begin{corollary}\label{fom}  Let $f$ be as in {\bf Theorem F}. If
for all $\chi+\eta+\zeta={\bf1},$ with $\chi\ne{\bf 0},$ all the
values of type {\rm (\ref{esremd})} are finite, $F\in L(\mathbb
R^d)$ if and only if

\begin{eqnarray*}\int\limits_{\mathbb{R}^d_+}|f(x)|\prod\limits_{j=1}^d
\frac{1}{x_j}\,dx<\infty.\end{eqnarray*}\end{corollary}

Theorem D and results from \cite{GM1} are immediate corollaries of
Theorem F. And, of course, Theorem E can easily be deduced, say,
from Corollary \ref{fom}.

\bigskip

\section{Proofs}\label{proofs}

We give, step by step, proofs of the main results formulated in
Introduction.

\subsection{Proof of Theorem \ref{th1}.}

To prove {\bf a)}, we apply the first part of Theorem C.

Denoting $h(p)=\pi2^{-p},$ $p\in\mathbb Z,$ and

\begin{eqnarray}\label{dl}\Delta(h)=\left(\int\limits_{\mathbb R}
|f(t+h)-f(t-h)|^2dt\right)^{1/2},                 \end{eqnarray}
we are going to prove that

\begin{eqnarray}\label{os2}\sum\limits_{p=0}^\infty 2^{p/2}\Delta(h(p))
+\sum\limits_{p=1}^\infty2^{-p/2}\Delta(h(-p))<\infty.\end{eqnarray}
It is obvious that for $h>0$

\begin{eqnarray}\label{pr1} |f(t+h)-f(t-h)|\le2f_0(\min|t\pm
h|), \end{eqnarray}

\begin{eqnarray}\label{pr3} \qquad|f(t+h)-f(t-h)|=|\int\limits_{t-h}^{t+h}
f'(s)\,ds|\le 2hf_1(\min|t\pm h|),               \end{eqnarray}
and for $|t|\le 3h$

\begin{eqnarray}\label{pr4} \qquad|f(t+h)-f(t-h)|\le\int\limits_{t-h}^{t+h}
f_1(s)\,ds|\le 2\int\limits_0^{4h}f_1(t)\,dt.    \end{eqnarray}

The proof will be divided into several steps.

{\it Step 1.} To separately study the behavior near the origin and
near infinity, we represent $f$ as the sum of two functions
$\varphi$ and $\psi$ with similar properties. First, let
$\varphi(t)=f(t)$ when $|t|\le2\pi,$ while for $|t|\ge2\pi$ it is
$f(t)(3-\frac{|t|} {\pi})_+.$ Consequently,
$\psi(t)=f(t)-\varphi(t).$

Monotone majorants of the absolute values of these functions and
their derivatives satisfy the inequalities

\begin{eqnarray}\label{in5} \varphi_0(t)\le f_0(t),\quad
\varphi_1(t)\le f_1(t)+\frac{1}{\pi}f_0(t),\end{eqnarray} and

\begin{eqnarray}\label{in6} \psi_0(t)\le f_0(t),\quad
\psi_1(t)\le\begin{cases} f_1(2\pi)+\frac{1}{\pi}f_0(2\pi), &
|t|\le3\pi,\\ f_1(t), & |t|\ge3\pi.\end{cases}\end{eqnarray}

{\it Step 2.} For the compactly supported function $\varphi$ the
second sum in (\ref{os2}) does not exceed

\begin{eqnarray*} 2\sum\limits_{p=1}^\infty 2^{-p/2}||\varphi||_2
\le2\sum\limits_{p=1}^\infty
2^{-p/2}\biggl(\int\limits_{-3\pi}^{3\pi}
|f(t)|^2dt\biggr)^{1/2}\le c_1 f_0(0).\end{eqnarray*} Further,

\begin{eqnarray}\label{in6.5} f_0(0)&=&\sup\limits_{t\in\mathbb R}
|f(t)|\le\sup\limits_{|t|\le2}|f(t)|+f_0(2)\nonumber\\
&\le&\sup\limits_{|t|\le2}|f(t)-f(2\mbox{\rm sign}t)|+2f_0(2)\nonumber\\
&\le&\int\limits_0^2 f_1(t)\,dt+2f_0(2)\le2\int\limits_0^1 f_1(t)\,dt+2f_0(2)\nonumber\\
&\le&\frac{2}{\ln2}\int\limits_0^1 f_1(t)\ln\frac{2}{t}\,dt+
\frac{2} {\ln2}\int\limits_1^2
\frac{f_0(t)}{t}\,dt\le\frac{2}{\ln2}(A_1+A_0).
\end{eqnarray}

{\it Step 3.} To estimate the first sum in (\ref{os2}), for
$\varphi$ equal to

\begin{eqnarray*}\sum\limits_{p=0}^\infty 2^{p/2}\biggl(
\int\limits_{|t|\le3\pi+h(p)}|\varphi(t+h(p))-\varphi(t-h(p)|^2dt
\biggr)^{1/2},                                   \end{eqnarray*}
we split the integral in (\ref{dl}) into the two ones: over
$|t|\le3h(p)$ and over $3h(p)\le|t|.$

{\it Step 3.1.} Using (\ref{pr4}), we obtain for the first part

\begin{eqnarray*}&\quad&\sum\limits_{p=0}^\infty 2^{p/2}\biggl(
\int\limits_{|t|\le3h(p)}|\varphi(t+h(p))-\varphi(t-h(p)|^2dt\biggr)^{1/2}\\
&\le&2\sum\limits_{p=0}^\infty 2^{p/2}\biggl(
\int\limits_{-3h(p)}^{3h(p)}\biggl(\int\limits_0^{4h(p)}\varphi_1(t)\,dt\biggr)^2
du\biggr)^{1/2}\\&=&2\sqrt{6\pi}\sum\limits_{p=0}^\infty
\int\limits_0^{4h(p)}\varphi_1(t)\,dt.\end{eqnarray*}

We will systematically need to pass from sums to integrals. In the
following simple inequalities one can pass to the limit as
$m\to\infty$ and $n\to\infty.$

If $g$ increases on $(0,\infty),$ and $n\ge m,$ then for any
$\alpha\in\mathbb R$

\begin{eqnarray}\label{in7}\sum\limits_{p=m}^n 2^{p\alpha}g(2^p)
&\le& \sum\limits_{p=m}^n \int\limits_{2^p}^{2^{p+1}}
t^{\alpha-1}g(t)\,dt\left(\int\limits_{2^p}^{2^{p+1}}
t^{\alpha-1}\,dt\right)^{-1}\nonumber \\
&=&\frac{\alpha}{2^\alpha-1}\int\limits_{2^m}^{2^{n+1}}
t^{\alpha-1}\varphi(t)\,dt.                     \end{eqnarray}

And if $g$ decreases on $(0,\infty),$ and $n\ge m,$ then for any
$\alpha\in\mathbb R$

\begin{eqnarray}\label{in8}\sum\limits_{p=m}^n 2^{p\alpha}g(2^p)&\le&
\sum\limits_{p=m}^n \int\limits^{2^p}_{2^{p-1}}
t^{\alpha-1}g(t)\,dt\left(\int\limits^{2^p}_{2^{p-1}}
t^{\alpha-1}\,dt\right)^{-1}\nonumber \\
&=&\frac{2^\alpha\alpha}{2^\alpha-1}\int\limits_{2^{m-1}}^{2^n}
t^{\alpha-1}g(t)\,dt.                       \end{eqnarray}

Passing to the integral, changing the order of integration and
substitution $s=8\pi u,$ we obtain, times a constant,

\begin{eqnarray}\label{in9}&\quad&\int\limits_0^{8\pi}\Big(\int\limits_0^{t}
\varphi_1(s)\,ds\Big)\,\frac{dt}{t}=\int\limits_0^{8\pi}\varphi_1(u)\,du
\int\limits_u^{8\pi}\frac{dt}{t}\nonumber\\
&=&\int\limits_0^{8\pi}\varphi_1(u)\ln\frac{8\pi}{u}\,du= 8\pi
\int\limits_0^1\varphi_1(8\pi u)\ln\frac{1}{u}\,du\le 8\pi A_1.
\end{eqnarray}

{\it Step 3.2.} Further, the second part of the sum does not
exceed, by (\ref{pr3}),

\begin{eqnarray*}&\quad&\sum\limits_{p=0}^\infty 2^{p/2}\left(
\int\limits_{3h(p)\le|t|\le4\pi}|\varphi(t+h(p))-\varphi(t-h(p))|^2
dt\right)^{1/2}\\
&\le&2\sum\limits_{p=0}^\infty 2^{p/2}h(p)\left(
\int\limits_{3h(p)\le|t|\le4\pi}\varphi_1^2(\min|t\pm h(p)|)\,dt
\right)^{1/2}.                                   \end{eqnarray*}
It follows from this, by taking into account the evenness of
$\varphi_1,$ the inequality (\ref{in7}) and the shift in the
integral $u=s-h(p),$ that

\begin{eqnarray*} &\quad&4\pi\sum\limits_{p=0}^\infty 2^{-p/2}\left(
\int\limits_{\pi2^{-p}}^{4\pi}\varphi_1^2(s)\,ds\right)^{1/2}\le
c_2\int\limits_0^{4\pi} t^{-1/2}\left(
\int\limits_u^{4\pi}\varphi_1^2(u)\,du\right)^{1/2}dt\\
&=&c_2\int\limits_0^{4\pi}\frac{1}{t^{1/2}\ln(8\pi/t)}\,\ln(8\pi/t)\left(
\int\limits_t^{4\pi}\varphi_1^2(u)\,du\right)^{1/2}dt.\end{eqnarray*}
Applying the Cauchy-Schwarz-Bunyakovskii inequality, we estimate
the above through

\begin{eqnarray*}&\quad&c_3\left(\int\limits_0^{4\pi}\ln^2\frac{8\pi}{t}
\int\limits_t^{4\pi}\varphi_1^2(u)\,du\,dt\right)^{1/2}\\
&=&c_3\left(\int\limits_0^{4\pi}\varphi_1^2(u)
\int\limits_0^u\ln^2\frac{8\pi}{t}\,dt\,du\right)^{1/2}\\
&\le&30c_3\left(\int\limits_0^{4\pi}\varphi_1^2(t)t\ln^2(8\pi/t)\,dt\right)^{1/2}\\
&\le&c_4\left(\int\limits_0^1\varphi_1^2(t)t\ln^2(2/t)\,dt\right)^{1/2}.\end{eqnarray*}
Here, as above while establishing (\ref{in9}), monotonicity of
$\varphi_1$ is used. Again by this, along with (\ref{in5}) and
(\ref{in6.5}), we have

\begin{eqnarray*}\varphi_1(t)\ln(2/t)\le2\int\limits_{t/2}^t
\varphi_1(u)\ln(2/u)\,du\le 2A_1+c_5f_0(0)\le c_6(A_0+A_1).
\end{eqnarray*} Therefore,

\begin{eqnarray*} &\left(\int\limits_0^1\varphi_1^2(t)t\ln^2(2/t)
\,dt\right)^{1/2}\le\sqrt{c_6(A_0+A_1)}\left(\int\limits_0^1\varphi_1(t)
\ln(2/t)\,dt\right)^{1/2}\\&\le\sqrt{c_6(A_0+A_1)}\sqrt{A_1+c_7f_0(0)}
\le c_{8}(A_0+A_1).\end{eqnarray*}

{\it Step 4.} Let us proceed to the function
$\psi(t)=f(t)-\varphi(t).$ It vanishes for $|t|\le2\pi$ and
coincides with $f$ for $|t|\ge3\pi.$ Indeed, for $|t|\le3\pi$ we
have $\psi(t)=f(t)(\frac{|t|}{\pi} -2)_+;$ see also (\ref{in6}).

{\it Step 4.1.} To prove the validity of (\ref{os2}), let us start
with the second sum. We split the integral in (\ref{dl}) into the
two ones: over $|t|\le3h(-p)$ and over $3h(-p)\le|t|.$ We have
(see (\ref{pr1}) and (\ref{pr3}))

\begin{eqnarray*}&\quad&\sum\limits_{p=1}^\infty 2^{-p/2}\,\biggl(
\,\int\limits_{3h(-p)\le|t|}|\psi(t+h(-p))-\psi(t-h(-p))|^2dt\biggr)^{1/2}\\
&\le&2\sum\limits_{p=1}^\infty2^{-p/2}\,\biggl(\,\int\limits_{3h(-p)\le|t|}h(-p)\\
&\times&\psi_0(\min|t\pm h(-p)|)\psi_1(\min|t\pm h(-p)|)\,dt\biggr)^{1/2}\\
&=&2\sqrt{2\pi}\sum\limits_{p=1}^\infty\,\biggl(\,\int\limits_{3h(-p)}^\infty
\psi_0(t-h(-p))\psi_1(t-h(-p))\,dt\biggr)^{1/2}.
\end{eqnarray*} The last integral is

\begin{eqnarray}\label{f4}&\quad&2\sqrt{2\pi}\sum\limits_{p=1}^\infty\,\biggl(\,
\int\limits_{2\pi2^p}^\infty \psi_0(t)\psi_1(t)\,dt\biggr)^{1/2}\nonumber\\
&\le&c_9 \int\limits_1^\infty\biggl(\, \int\limits_u^\infty
\psi_0(t)\psi_1(t)\,dt\biggr)^{1/2}\frac{du}{u}=c_9A_{0,1}.\end{eqnarray}

{\it Step 4.2.} The rest of the second sum is

\begin{eqnarray*}&\quad&\sum\limits_{p=1}^\infty 2^{-p/2}\,\biggl(
\,\int\limits_{3h(-p)\ge|t|}|\psi(t+h(-p))-\psi(t-h(-p))|^2dt\biggr)^{1/2}\\
&\le&\sum\limits_{p=1}^\infty 2^{-p/2}\,\biggl[\biggl(
\,\int\limits_{3h(-p)\ge|t|}|\psi(t+h(-p))|^2dt\biggr)^{1/2}\\&+&\biggl(
\,\int\limits_{3h(-p)\ge|t|}|\psi(t-h(-p))|^2dt\biggr)^{1/2}\biggr]\\
&\le&2\sum\limits_{p=1}^\infty 2^{-p/2}\,\biggl(
\,\int\limits_{4h(-p)\ge|t|}|\psi(t)|^2dt\biggr)^{1/2}.
\end{eqnarray*}
The right-hand side does not exceed, times a constant,

\begin{eqnarray}\label{f5}\qquad\int\limits_{2\pi}^\infty\biggl(\, \int\limits_{|t|\le u}
|\psi(t)|^2dt\biggr)^{1/2}\frac{du}{u^{3/2}}\le2\int\limits_{2\pi}^\infty\biggl(\,
\int\limits_{2\pi}^u \psi_0(t)^2dt\biggr)^{1/2}\frac{du}{u^{3/2}}.
\end{eqnarray}
We can consider in the sequel $\psi_0(2\pi)>0,$ since otherwise,
if $\psi_0(2\pi)=0$ (or, equivalently, if $\psi_1(2\pi)=0$) we
have $\psi(t)\equiv0.$

Integrating by parts in the last integral, we obtain

\begin{eqnarray*} &\quad&\biggl[-\frac{2}{\sqrt t}\biggl(\int\limits_{2\pi}^t
\psi_0^2(s)\,ds\biggr)^{1/2}\biggr]_{2\pi}^\infty\\
&+&\int\limits_{2\pi}^\infty\biggl(\int\limits_{2\pi}^t
\psi_0^2(s)\,ds\biggr)^{-1/2}\psi_0^2(t)\frac{dt}{\sqrt{t}}.\end{eqnarray*}

Integrated terms vanish at infinity by L'Hospital rule, say. For
$t>2\pi,$ by monotonicity of $\psi_0,$

\begin{eqnarray*}\biggl(\int\limits_{2\pi}^u \psi_0^2(s)\,ds\biggr)^{-1/2}
\le\frac{1}{\psi_0(u)\sqrt{u-2\pi}}.             \end{eqnarray*}
Since $u-2\pi\ge u/3$ when $u\ge3\pi,$ we arrive at the upper
bound

\begin{eqnarray*}\int\limits_{2\pi}^\infty\frac{\psi_0(u)}{\sqrt{u(u-2\pi)}}\,du
\le\psi_0(2\pi)\int\limits_{2\pi}^{3\pi}\frac{du}{\sqrt{u(u-2\pi)}}+
\sqrt 3 \int\limits_{3\pi}^\infty\frac{\psi_0(u)}{u}\,du.
\end{eqnarray*} It remains to take into account that

\begin{eqnarray*} \psi_0(2\pi)\le2\int\limits_{\pi}^{2\pi}\frac{\psi_0(t)}{t}\,dt
\le 2A_0,\end{eqnarray*} since also $\psi_0(t)\le f_0(t).$

{\it Step 4.3.} Let us go on to the first sum in (\ref{os2}).
Since $f(t)=0$ for $|t|\le2\pi,$ it is equal, with
$h(p)=\pi2^{-p}\in(0,\pi],$ to

\begin{eqnarray}\label{re6}\quad&\quad&\sum\limits_{p=0}^\infty 2^{p/2}\left(
\int\limits_{\max|t\pm h(p)|\ge2\pi}|\psi(t+h(p))-\psi(t-h(p))|^2dt\right)^{1/2}\nonumber\\
\qquad&=&\sum\limits_{p=0}^\infty2^{p/2}\left(\int\limits_{|t|\ge2\pi+
h(p)}|\psi(t+h(p))-\psi(t-h(p))|^2dt\right)^{1/2}.\end{eqnarray}

Here we again split the integral into two, but in a different way:
over $|t|\le 10\sqrt{2^p}$ and over $|t|\ge 10\sqrt{2^p}.$ In the
first one, we apply (\ref{pr3}):

\begin{eqnarray*}&\quad&\sum\limits_{p=0}^\infty 2^{p/2}\biggl(
\int\limits_{2\pi+h(p)\le|t|\le 10\sqrt{2^p}} \psi_1^2(\min|t\pm
h(p)|)\,dt \biggr)^{1/2}2h(p)\\&\le& c_{10}\sum\limits_{p=0}^\infty
2^{-p/2}\biggl( \int\limits_{2\pi}^{2\sqrt{\pi2^p}}\psi_1^2(t)\,dt
\biggr)^{1/2}.\end{eqnarray*}

Passing to the integral, we get (cf. (\ref{f5}))

\begin{eqnarray*}\int\limits_{2\pi}^\infty\biggl(\,\int\limits_{2\pi}^{\sqrt{u}}
\psi_1(t)^2dt\biggr)^{1/2}\frac{du}{u^{3/2}}.\end{eqnarray*}
Substituting $\sqrt{u}\to u,$ we have

\begin{eqnarray*}\int\limits_{2\pi}^\infty\biggl(\,
\int\limits_{2\pi}^u \psi_1(t)^2dt\biggr)^{1/2}\frac{du}{u^2}.
\end{eqnarray*}

Repeating similar estimations as for $\psi_0$ (see (\ref{f5}) and
further), we obtain, times a constant,

\begin{eqnarray*}\psi_1(2\pi)+\int\limits_{3\pi}^\infty\frac{\psi_1(t)}
{t^{3/2}}\,dt\le c_{11}\psi_1(2\pi)\le c_{11}(f_1(2\pi)+f_0(2\pi)).
\end{eqnarray*}
By monotonicity of $f_0$ and $f_1$ the last bound is obviously
controlled by $A_0+A_1$ (see (\ref{in6.5})).

{\it Step 4.4.} In the remained sum all is similar to getting
(\ref{f4}). The only difference is that we get

\begin{eqnarray*}\int\limits_1^\infty\biggl(\, \int\limits_{\sqrt{u}}^\infty
\psi_0(t)\psi_1(t)\,dt\biggr)^{1/2} \frac{du}{u} \end{eqnarray*}
as an upper bound. But after substituting $\sqrt{u}\to u,$ we
obtain exactly $A_{0,1},$ which completes the proof of {\bf a)}.

The proof of {\bf b)} is also based on the first part of Theorem
C. We first replace (\ref{pr3}) with

\begin{eqnarray}\label{pr12} \qquad|f(t+h)-f(t-h)|=|\int\limits_{t-h}^{t+h}
f'(s)\,ds|\le 2hf_\infty(\max|t\pm h|).            \end{eqnarray}
It follows from $A_\d<\infty$ that for $|t|\ge2\pi$ and for the
same $\delta\in(0,1)$

\begin{eqnarray*}|f(t)|\le|f_0(t)|\le\left(\frac{A_\d^{1+\delta}}{f_\infty(4\pi)}
\right)^{1/\delta}\,|t|^{-1/\delta},      \end{eqnarray*}
therefore

\begin{eqnarray*} ||f||_2^2=\int\limits_{|t|\ge2\pi}|f(t)|^2dt\le
c_{12}\left(\frac{A_\d^{1+\delta}}{f_\infty(4\pi)}\right)^{2/\delta}.
\end{eqnarray*}
Hence the second sum in (\ref{os2}) does not exceed

\begin{eqnarray}\label{in13}2||f||_2\sum\limits_{p=1}^\infty2^{-p/2}\le
c_{13}\left(\frac{A_\d^{1+\delta}}{f_\infty(4\pi)}\right)^{1/\delta}.
\end{eqnarray}

Let $\delta_1=\frac{1-\delta}{1+\delta}.$ Clearly,
$\delta_1\in(0,1).$ Applying to the first sum in (\ref{os2}) (cf.
also (\ref{re6})) simultaneously (\ref{pr1}) and (\ref{pr12}), we
bound it with

\begin{eqnarray*}&\quad&2\sum\limits_{p=0}^\infty 2^{p/2}\biggl(
\int\limits_{|t|\ge2\pi+h(p)} f_0^{1-\delta_1}(\min|t\pm h(p)|)\\
&\times& h(p)^{1+\delta_1}f_\infty^{1+\delta_1}(\max|t\pm h(p)|)\,
dt\biggr)^{1/2}\\
&=&4\sum\limits_{p=0}^\infty 2^{p/2}h(p)^{(1+\delta_1)/2}\left(
\int\limits_{2\pi}^\infty f_0^{1-\delta_1}(t)f_\infty^{1+\delta_1}
(t+2h(p))\,dt\right)^{1/2}\\
&\le&\gamma_1(\delta)\left( \int\limits_{2\pi}^\infty
f_0^{1-\delta_1}(t)f_\infty^{1+\delta_1}(t+2\pi)\,dt\right)^{1/2}.\end{eqnarray*}
The choice of $\delta_1$ and assumptions of the theorem yield for
$t\ge2\pi$

\begin{eqnarray*}f_0^{1-\delta_1}(t)f_\infty^{1+\delta_1}(t+2\pi)
=(f_0^{\delta}(t)f_\infty(t+2\pi))^{2/(1+\delta)}\le \frac{A_\d^2}
{t^{\frac{2}{1+\delta}}},                       \end{eqnarray*}
and the first sum does not exceed $\gamma_2(\delta)A_\d.$
Combining it with (\ref{in13}) gives the desired estimate.

It is also worth mentioning that
$\max\limits_{2\pi\le|t|\le4\pi}|f(t)|\le 2\pi f_\infty(4\pi).$

The proof of the theorem is complete.   \hfill$\Box$

\subsection{Proof of Corollary \ref{sl}.}

If $\a>0$ and $\b\ge0 $ with $\a+\b>1,$ we apply {\bf a)} of
Theorem \ref{th1}, while if $\b<0$ but still $\a+\b>1,$ the
assertion {\bf b)} of Theorem \ref{th1} is applicable (one can
take any $\d$ satisfying $\d\a+\b=1$).

There is a counterexample in the multiple case as well
(\cite[7.4]{Stein}). For $d=1,$ the function

\begin{eqnarray*}
g(x)=\frac {e^{i|t|^{\a_1}}}{(1+|t|^2)^{\b_1}}  \end{eqnarray*}
does not belong to $A(\BR)$  if $\a_1\ne 1$ and $4\b_1<\a_1.$

If $\a\ne \b$ we set $\b_1=\frac \a2 $ and $\a_1=\a-\b+1.$ Then
$4\b_1-\a_1=\a+\b-1<0$ and for $|t|\to \infty$

\begin{eqnarray*}  g(t)=O\Big(\frac {1}{|t|^\a}\Big),\quad
g^{\prime}(t)=O\Big(\frac{1}{|t|^\b)}\Big).      \end{eqnarray*}
If $\a=\b$ and $\a+\b<1,$ one can a larger $\a$ so that the sum to
still be less than 1 and then to make use of the above argument.

\subsection{Proof of Theorem \ref{newmu}.}

The proof will go along the same lines, or, more precisely, the same
steps, as that of {\bf a)} of Theorem \ref{th1} does. We first
represent the given function $f$ as the sum of two functions
$\varphi$ and $\psi$ so that $\varphi$ to be compactly supported and
near the origin coincide with $f,$ while $\psi$ correspondingly
vanishes near the origin and coincides with $f$ near infinity. Thus,
let $\varphi(x)=f(x)$ when $|x_j|\le2\pi,$ $j=1,2,...,d.$ Further,
when $|x_j|\ge2\pi$ and $|x_j|\ge |x_k|$ for all $k=1,2,...,d,$
$\varphi(x)=f(x)(3-\frac{|x_j|}{\pi})_+,$ and this is for each
$j=1,2,...,d.$ Correspondingly, $\psi(x)=f(x)-\varphi(x).$

Extensions of (\ref{pr1}), (\ref{pr3}) and (\ref{pr4}) are as
follows. Let $r=\chi_1+...+\chi_d,$ $h=(h_1,...,h_d),$ and
$x_{\chi,\pm h}$ be the vector $x$ with $x_j$ replaced by $|x_j\pm
h_j|$ for $j:\chi_j=1.$

\begin{eqnarray}\label{pr1m} |\prod\limits_{j:\chi_j=1}\Delta_{h_j}f(x)|
\le2^rf_{\chi,{\bf 0}}(\min x_{\chi,\pm h}), \end{eqnarray}

\begin{eqnarray}\label{pr3m}|\prod\limits_{j:\chi_j=1}\Delta_{h_j}f(x)|
\le2^rf_{{\bf 0},\chi}(\min x_{\chi,\pm h}),
\end{eqnarray} and when $|x_j|\le 3h_j,$ $j:\chi_j=1$

\begin{eqnarray}\label{pr4m} |\prod\limits_{j:\chi_j=1}\Delta_{h_j}f(x)|
\le 2^r\int\limits_{\prod[0,4h_j], \atop j:\chi_j=1}f_{{\bf
0},\chi}(x)\,dx_\chi.  \end{eqnarray}

Now, the result will follow from the next two propositions, in
which $\varphi$ and $\psi$ are treated separately. The first one
corresponds to the case in (\ref{muco}) when $\chi=\eta={\bf 0}$
and $\zeta={\bf 1}.$

\begin{proposition}\label{br} Let $\varphi(x)\in C_0(\mathbb R^d)$
be supported on $\{x: |x_j|\le3\pi,j=1,...,d\}.$ Let $\varphi$ and
its partial derivatives $D^\eta \varphi,$ ${\bf 0}\le\eta<{\bf 1},$
be locally absolutely continuous on $(\mathbb R\setminus \{0\})^d$
in each variable. If

\begin{eqnarray}\label{lnc} \int\limits_0^1\ln\frac{2}{x_1}\dots
\int\limits_0^1\ln\frac{2}{x_d} |\varphi_1(x)|\,dx<\infty,
\end{eqnarray} then $\varphi\in A(\mathbb R^d).$ \end{proposition}

\begin{proof} We again use a) of Theorem C. By this, we have to
estimate the sums of the form

\begin{eqnarray}\label{twos} \sum\limits_{1\le s_i<\infty,\atop i: \chi_i=0}
2^{-\frac 12\sum\limits_{i: \chi_i=0} s_i}\sum\limits_{0\le
s_j<\infty,\atop j: \chi_j=1} 2^{\frac12\sum
\limits_{j:\chi_j=1}s_j}||\D_{H(s)}(\varphi)||_2<\infty,\end{eqnarray}
where $H(s)$ is the $d$-dimensional vector with the entries
${\pi}{2^{s_i}}$ when $\chi_i=0$ and ${\pi}{2^{-s_j}}$ when
$\chi_j=1.$ Proceeding to the first sum, we have

\begin{eqnarray*} &\quad&\sum\limits_{1\le s_i<\infty,\atop i:\chi_i=0}
2^{-\frac 12\sum\limits_{i:\chi_i=0}s_i}||\D_{H(s)}(\varphi)||_2\\
&\le&2^d\sum\limits_{1\le s_i<\infty,\atop i: \chi_i=0} 2^{-\frac
12\sum\limits_{i:\chi_i=0} s_i}||
\prod\limits_{j:\chi_j=1}\D_{\pi2^{-s_j}}\varphi||_2\\
&\le&2^d\sum\limits_{1\le s_i<\infty,\atop i:\chi_i=0}
2^{-\frac12\sum_{i:\chi_i=0}s_i} \biggl(\int\limits_{|x_i|
\le3\pi,\atop i:\chi_i=0}\int\limits_{|x_i|\le3\pi, \atop
i:\chi_i=1}|\prod\limits_{j:\chi_j=1}\D_{\pi2^{-s_j}}f(x)|^2dx\biggr)^{1/2}\\
&\le&c_{14}\biggl(\int\limits_{|x_i|\le3\pi,\atop i: \chi_i=1}
\biggl( \prod\limits_{j:\chi_j=1}\D_{\pi2^{-s_j}}f\biggr)_{{\bf
1}-\chi,{\bf 0}}(x_\chi^0)^2 dx_\chi\biggr)^{1/2},\end{eqnarray*}
where $x_\chi^0$ is $x$ with zero entries in place of $x_i$ when
$\chi_i=0.$

Going on to the second sum in (\ref{twos}), we have to estimate

\begin{eqnarray*}\sum\limits_{0\le s_j<\infty,\atop j:\chi_j=1} 2^
{\frac12\sum\limits_{j:\chi_j=1}s_j}\left(\int\limits_{|x_i|\le3\pi,\atop
i: \chi_i=1} \biggl( \prod\limits_{j:\chi_j=1}\D_{\pi2^{-s_j}}
f\biggr)_{{\bf 1}-\chi,{\bf 0}}(x_\chi^0)^2 dx_\chi\right)^{1/2}.
\end{eqnarray*}
But this is estimated similarly to the one-dimensional case, with
calculations repeated in each $j$-th variable for $j:\chi_j=1.$
This completes the proof. \hfill\end{proof}

Let us go on to $\psi.$

\begin{proposition}\label{ni} Let $\psi(x)\in C_0(\mathbb R^d)$
vanish on $\{x: |x_j|\le2\pi,j=1,...,d\}.$ Let $\psi$ and its
partial derivatives $D^\eta \psi,$ ${\bf 0}\le\eta<{\bf 1},$ be
locally absolutely continuous on $(\mathbb R\setminus \{0\})^d$ in
each variable. Let also partial derivatives $D^\eta \psi,$ ${\bf
0}<\eta\le{\bf 1}$ be almost everywhere bounded out of any
neighborhood of each coordinate hyperplane. If

\begin{eqnarray*}\biggl\|\int\limits_1^\infty...\int\limits_1^\infty
\biggl(\int\limits_{\prod\limits_{j:\eta_j=1}[u_j,\infty)}f_{\chi+\eta,\zeta}(x)
f_{\chi,\eta+\zeta}(x)\,dx_\eta\biggr)^{1/2}\prod\limits_{i:\chi_i=1\atop
{\rm or}\ \eta_i=1}\frac {du_i}{u_i}\biggr\|_{L_\infty(\mathbb
R_\zeta)} <\infty                                 \end{eqnarray*}
for all $\chi,$ $\eta$ and $\zeta$ such that $\chi+\eta+\zeta={\bf
1}=\{1,1,...,1\},$ then $f\in A(\mathbb R^d).$\end{proposition}

\begin{proof} We again use {\bf a)} of Theorem C. The steps of the proof
will be similar to the four sub-steps of Step 4 in the
one-dimensional proof. The first two are concerned with the sums
where $1\le s_j<\infty$ and the factor $2^{s_j/2}$ stays before
the $L^2$ norm, while in the next two $0\le s_j<\infty$ and the
factor $2^{-s_j/2}$ stays before the $L^2$ norm. Each two steps
correspond to splitting the relevant integrals over
$|x_j|\le3h(-s_j)$ and $3h(-s_j)\le|x_j|$ and
$|x_j|\le\sqrt{3h(-s_j)}$ and $\sqrt{3h(-s_j)}\le|x_j|,$
respectively. For each case separately the estimates are similar,
we just repeat the same calculations in each variable involved.
Let us indicate certain points in this procedure.

First, (\ref{pr1m}), (\ref{pr3m}) and (\ref{pr4m}) are used
instead of (\ref{pr1}), (\ref{pr3}) and (\ref{pr4}), respectively.

Then, in the one-dimensional case the order of estimations was not
important, since they were completely independent. Here all types
of estimates can be applicable at once, each with respect to
certain group of variables assigned by $\chi,$ $\eta$ and $\zeta.$
We start with an analog of Step 4.2. Let us show how to deal with
integrated terms while integrating by parts. Without loss of
generality, we can consider

\begin{eqnarray}\label{mod} \int\limits_{2\pi}^\infty\int\limits_{2\pi}^\infty
\left(\int\limits_{2\pi}^x\int\limits_{2\pi}^y
F^2(s,t)\,ds\,dt\right)^{1/2}\frac{dx}{x^{3/2}}\,\frac{dy}{y^{3/2}}
\end{eqnarray}
to be a model case, with $F$ bounded and vanishing at infinity.
Integrating by parts in $y,$ we obtain

\begin{eqnarray*}&\quad&\biggl[-\frac{2}{\sqrt
y}\biggl(\int\limits_{2\pi}^x\int\limits_{2\pi}^y
F^2(s,t)\,ds\,dt\biggr)^{1/2}\biggr]_{2\pi}^\infty\\
&+&\int\limits_{2\pi}^\infty\biggl(\int\limits_{2\pi}^x\int\limits_{2\pi}^y
F^2(s,t)\,ds\,dt\biggr)^{-1/2}\biggl(\int\limits_{2\pi}^x
F^2(s,y)\,ds\,dt\biggr)\frac{dy}{y^{1/2}}.    \end{eqnarray*}
Since $F$ is bounded, we can pass to the limit under the integral
sign while using the L'Hospital rule as above. The estimates $y$
are exactly the same as in dimension one, and then we just repeat
these in $x.$ By this, (\ref{mod}) is controlled by

\begin{eqnarray*}F(2\pi,2\pi)+\int\limits_{2\pi}^\infty\frac{F(x,2\pi)}{x}\,dx+
\int\limits_{2\pi}^\infty\frac{F(2\pi,y)}{y}\,dy+
\int\limits_{2\pi}^\infty\int\limits_{2\pi}^\infty
\frac{F(x,y)}{xy}\,dx\,dy.\end{eqnarray*}

Further, we fulfil an analog of Step 4.3 for the corresponding
group of variables. And after that we make use of analogs of Steps
4.1 and 4.4 that are, in essence, the same. We finally arrive at
the desired estimate with $A_{\chi,\eta,\zeta}.$ Of course, for
some choices of $\chi,$ $\eta$ and $\zeta$ not all types of
estimates to be involved; there are extremal cases when only one
type of estimates, that is, all sums and all integrals are of one
type for all the variables, is fulfilled.

The interplay of the steps is the main problem in the multivariate
extension, in dimension one all the steps are independent. To make
such an interplay more transparent, let us consider all the
details with easier notation in dimension two. More precisely, we
consider the cases where first Steps 4.1 and 4.2 and then Steps
4.1 and 4.4 occur simultaneously. Since both cases are
two-dimensional we switch the notation from that with subscripts
to the one with different letters.

Thus, let us estimate

\begin{eqnarray*} &\quad&\sum\limits_{p=1}^\infty 2^{-p/2}
\sum\limits_{q=1}^\infty 2^{-q/2}\biggl(\ \int\limits_{|s|
\ge3\pi2^p}\int\limits_{|t|\le3\pi2^q}|f(s+h(-p),t+h(-q))\\
&-&f(s-h(-p),t+h(-q))-f(s+h(-p),t-h(-q))\\&+&f(s-h(-p),t-h(-q))|^2
ds\,dt\biggr)^{1/2}.\end{eqnarray*}

Denoting, for brevity, $\Omega(s,t)=f(s+h(-p),t)f(s-h(-p),t),$ we
first estimate, more or less along the same lines as in Step 4.2,

\begin{eqnarray*}\sum\limits_{q=1}^\infty 2^{-q/2}\biggl(\
\int\limits_{|t|\le3\pi2^q}\int\limits_{|s|\ge3\pi2^p}
|\Omega(s,t+h(-q))-\Omega(s,t-h(-q))|^2ds\,dt\biggr)^{1/2}\end{eqnarray*}

\begin{eqnarray*}&\le&\sum\limits_{q=1}^\infty 2^{-q/2}\biggl(\
\int\limits_{|t|\le3\pi2^q}\int\limits_{|s|\ge3\pi2^p}
|\Omega(s,t+h(-q))|^2ds\,dt\biggr)^{1/2}\\
&+&\sum\limits_{q=1}^\infty 2^{-q/2}\biggl(\
\int\limits_{|t|\le3\pi2^q}\int\limits_{|s|\ge3\pi2^p}
|\Omega(s,t-h(-q))|^2ds\,dt\biggr)^{1/2}\\
&\le&2\sum\limits_{q=1}^\infty 2^{-q/2}\biggl(\
\int\limits_{|t|\le4\pi2^q}\int\limits_{|s|\ge3\pi2^p}
|\Omega(s,t)|^2ds\,dt\biggr)^{1/2}\\
&\le&c_{15}\int\limits_{2\pi}^\infty\biggl(\int\limits_{|t|\le y}
\int\limits_{|s|\ge3\pi2^p}|\Omega(s,t)|^2ds\,dt\biggr)^{1/2}dy\\
&\le&c_{16}\int\limits_{2\pi}^\infty\biggl(\int\limits_{2\pi}^y
\int\limits_{|s|\ge3\pi2^p}|\Omega_0(s,t)|^2ds\,dt\biggr)^{1/2}
\frac{dy}{y^{3/2}}.                            \end{eqnarray*} Here
$\Omega_0$ means, for a moment, the same as $f_0$ in Theorem
\ref{th1} but with respect to one of the variables, $t.$ Similarly,
$\Omega_1$ is an analog of $f_1.$ Exactly as in the Step 4.2 we
bound the right-hand side by

\begin{eqnarray}\label{ts}&\quad& c_{17}\biggl(\int\limits_0^1\biggl(
\int\limits_{|s|\ge3\pi2^p}|\Omega_1(s,y)|^2ds\biggr)^{1/2}\ln
\frac{2}{y}\,dy\nonumber\\
&+&\int\limits_1^\infty\biggl(\int\limits_{|s|\ge3\pi2^p}
|\Omega_0(s,y)|^2ds\biggr)^{1/2}\frac{dy}{y}\biggr).\end{eqnarray}
What remains is to estimate

\begin{eqnarray*}\sum\limits_{p=1}^\infty 2^{-p/2}\biggl(
\int\limits_{|s|\ge3\pi2^p}|F(s+h(-p),y)-F(s-h(-p),y)|^2ds\biggr)^{1/2},
\end{eqnarray*}
where $F$ denotes either $\Omega_0$ or $\Omega_1.$ But this is
exactly Step 4.1 above. Fulfilling it, we arrive at one of the
relations of type (\ref{muco}), with

\begin{eqnarray*}f_{(1,0),(0,1)}f_{(0,0),(1,1)}\end{eqnarray*}
inside while estimating the first summand in (\ref{ts}) and

\begin{eqnarray*}f_{(1,1),(0,0)} f_{(0,1),(1,0)}\end{eqnarray*}
for the second one.

Let us go on to to the combination of Steps 4.1 and 4.4. We are
proceeding to

\begin{eqnarray*} &\quad&\sum\limits_{p=1}^\infty 2^{-p/2}
\sum\limits_{q=0}^\infty 2^{q/2}\biggl(\
\int\limits_{|s|\ge3\pi2^p}\int\limits_{|t|\ge\sqrt{3\pi2^q}}
|f(s+h(-p),t+h(-q))\\
&-&f(s-h(-p),t+h(-q))-f(s+h(-p),t-h(-q))\\&+&f(s-h(-p),t-h(-q))|^2
ds\,dt\biggr)^{1/2}.                           \end{eqnarray*}
With the same notation $\Omega$ in hand and using Step 4.4, we
estimate

\begin{eqnarray*}\sum\limits_{q=0}^\infty 2^{-q/2}\biggl(\
\int\limits_{|t|\ge\sqrt{3\pi2^q}}\int\limits_{|s|\ge3\pi2^p}
|\Omega(s,t+h(-q))-\Omega(s,t-h(-q))|^2ds\,dt\biggr)^{1/2}\end{eqnarray*}

\begin{eqnarray*}&\le&c_{18}\sum\limits_{q=0}^\infty\biggl(\int\limits_{|s|\ge3\pi2^p}
\int\limits_{\sqrt{3\pi2^q}}^\infty\Omega_0(s,t)\Omega_1(s,t)\,dt
\,ds\biggr)^{1/2}\\
&\le&c_{19}\int\limits_{2\pi}^\infty\biggl(\int\limits_{|s|\ge3\pi2^p}
\int\limits_{y}^\infty\Omega_0(s,t)\Omega_1(s,t)\,dt
\,ds\biggr)^{1/2}\frac{dy}{y}.\end{eqnarray*}

Similarly,

\begin{eqnarray*} &\quad&\sum\limits_{p=1}^\infty 2^{-p/2}\biggl(
\int\limits_y^\infty\int\limits_{\sqrt{3\pi2^q}}^\infty\Omega_0(s,t)
\Omega_1(s,t)\,ds\,dt\biggr)^{1/2}\\
&\le&c_{20}\int\limits_{2\pi}^\infty\biggl(\int\limits_x^\infty
\int\limits_y^\infty f_{(1,1),(0,0)}(s,t)f_{(0,0),(1,1)}(s,t)\,ds
\,dt\biggr)^{1/2}\frac{dx}{x}\frac{dy}{y}.     \end{eqnarray*} This
is the needed bound. We just note that in the same way one could
have in the inner integral $f_{(0,1),(1,0)}f_{(1,0),(0,1)}$ or even

\begin{eqnarray*}\sqrt{f_{(1,1),(0,0)}f_{(0,1),(1,0)}
f_{(1,0),(0,1)}f_{(0,0),(1,1)}}.\end{eqnarray*}

The proof is complete.   \hfill\end{proof}

{\bf Acknowledgements.}

The authors gratefully acknowledge the support of the Gelbart
Research Institute for Mathematical Sciences at Bar-Ilan
University.

The work of the second author was also supported by the Ukranian
Fund for Fundamental Research Ukraine, Project F25.1/055.

We also wish to thank T. Shervashidze and his collaborators for
providing us with the text of almost inaccessible paper
\cite{Gab}.

\end{document}